\UseRawInputEncoding

\documentclass[12pt,reqno]{amsart}
\usepackage{amssymb,amsmath,graphicx,
	amsfonts,euscript}
\usepackage{color}
\usepackage{mathrsfs}

\let\f=\frac

\let\Om=\Omega

\def\R{\Bbb R}

\def\no{\noindent}

\def\endproof{\hphantom{MM}\hfill\llap{$\square$}\goodbreak}

\newcommand{\beq}{\begin{equation}}
\newcommand{\eeq}{\end{equation}}
\newcommand{\ben}{\begin{eqnarray}}
\newcommand{\een}{\end{eqnarray}}
\newcommand{\beno}{\begin{eqnarray*}}
\newcommand{\eeno}{\end{eqnarray*}}

\newtheorem{Theorem}{Theorem}[section]
\newtheorem{Definition}[Theorem]{Definition}
\newtheorem{Proposition}[Theorem]{Proposition}
\newtheorem{Lemma}[Theorem]{Lemma}

\newtheorem{Remark}[Theorem]{Remark}

\begin{document}

\title[interior and boundary regularity criteria of Navier-Stokes]{\bf Interior and Boundary Regularity Criteria for the 6D steady Navier-Stokes Equations}

\author{Shuai Li, Wendong Wang}

\address[Shuai, Li]{School of Mathematical Sciences, Dalian University of Technology, Dalian, 116024, P. R. China.} \email{leeshy@mail.dlut.edu.cn}

\address[Wendong, Wang]{School of Mathematical Sciences, Dalian University of Technology, Dalian, 116024, P. R. China.} \email{wendong@dlut.edu.cn}

%

\maketitle

\begin{abstract}
It is shown in this paper that suitable weak solutions to the 6D steady incompressible Navier-Stokes are H\"{o}lder continuous at $0$ provided that $\int_{B_1}|u(x)|^3dx+\int_{B_1}|f(x)|^qdx$ or $\int_{B_1}|\nabla u(x)|^2dx$+$\int_{B_1}|\nabla u(x)|^2dx\left(\int_{B_1}|u(x)|dx\right)^2+\int_{B_1}|f(x)|^qdx$ with $q>3$ is sufficiently small, which implies that the 2D Hausdorff measure of the set of singular points is zero.
For the boundary case,  we obtain that $0$ is regular provided that
$\int_{B_1^+} |u(x)|^3 dx +  \int_{B_1^+} |f(x)|^3 dx$ or $\int_{B_1^+} |\nabla u(x)|^2 dx +  \int_{B_1^+} |f(x)|^3 dx$ is sufficiently small.
These results improve previous  regularity theorems by Dong-Strain (\cite{DS}, Indiana Univ. Math. J., 2012), Dong-Gu (\cite{DG2}, J. Funct. Anal., 2014),  and Liu-Wang (\cite{LW}, J. Differential Equations, 2018), where either the smallness of the pressure or the smallness on all balls is necessary.
\end{abstract}



\setcounter{equation}{0}
\section{Introduction}

Consider the following 6D steady incompressible Navier-Stokes equations on $\Omega \subset \R^6$ as follows:
\begin{align} \label{eq:SNS}(\rm SNS)\,\, \left\{
\begin{aligned}
&-\Delta u+u\cdot \nabla u=-\nabla{\pi}+f,\\
&\quad\nabla\cdot u=0,\\
\end{aligned}
\right. \end{align}
where $u$ represents the fluid velocity field, $\pi$ is a scalar pressure.

The  $\varepsilon-$regularity analysis of the above equations is started by Struwe's question in \cite{St1,St2}, where he obtained partial regularity for $N=5$ by regularity methods of elliptic systems (c.f. Morrey \cite{Mo} and Giaqinta \cite{Gi}) and asked if analogous partial regularity results hold in spacial dimension $N>5.$
Later, the result of Struwe was extended to the boundary case by Kang \cite{Ka}.  Recently interior regularity results in 6D are obtained by Dong-Strain \cite{DS}, and
they proved $0$ is regular if
\beno
\limsup_{r\rightarrow0}r^{-2}\int_{B_r}|\nabla u|^2dx\leq \varepsilon_0.
\eeno
Moreover, similar boundary regularity results are obtained in Dong-Gu \cite{DG2} and Liu-Wang \cite{LW} by different methods, respectively.
For more developments, in a series of papers by Frehse and Ruzicka \cite{FR1,FR2,FR3,FR4}, the existence on a class of special regular solutions of (\ref{eq:SNS}) was obtained for the five-dimensional and higher dimensional case.  Gerhardt \cite{Ge} obtained the regularity of weak solutions under the  four-dimensional case.  More references, we refer to Li-Yang \cite{LY}  for the existence of regular solutions of high dimensional Navier-Stokes equations.   At last, we refer to  \cite{FS} by Farwig-Sohr  for existence and regularity criteria for weak solutions to inhomogeneous  Navier-Stokes equations.


Recall that these so-called $\varepsilon-$regularity criteria can be traced back to
the well-known work by Caffarelli-Kohn-Nirenberg \cite{CKN} for the analysis of suitable weak solutions of the three dimensional time-dependent Navier-Stokes equations, where they  showed that the set $\mathcal{S}$ of possible interior singular points of a suitable weak solution is one-dimensional parabolic Hausdorff measure zero by improving Scheffer's results in \cite{SV1,SV2,SV4}. More references on simplified proofs and improvements, we refer to Lin \cite{Lin}, Ladyzhenskaya-Seregin \cite{LS}, Tian-Xin \cite{TX}, Seregin \cite{Se}, Gustafson-Kang-Tsai \cite{GKT}, Vasseur \cite{Va}, Kukavica \cite{Ku}, Wang-Zhang \cite{WZ2} and the references therein.
Motivated by the recent interior regularity by Wolf \cite{Wolf}, where the author proved $\int_{Q_1}|u(x)|^3dx\leq \varepsilon_0$ in one scale can imply the regularity via pressure decomposition of Stokes equation. Also, we refer to Chae-Wolf \cite{CW} and \cite{JWZ,WWZ} for some recent progress. One can ask naturally:\\
{\bf  ``Whether the smallness of the velocity in a ball can ensure the interior or boundary regularity of the 6D  steady  Navier-Stokes equations?"}\\
In this note, we try to investigate this issue and answer these questions.

After  finishing  this  paper,  the  authors  have  become  to  know  that,  very  recently,  Cui \cite{Cui} showed that local interior regularity
and boundary regularity in one scale for the 5D steady Navier-Stokes equations via Campanato¡¯s method as Dong-Wang \cite{DW}. However, we considered the 6D case, which is the largest dimension, and used the Wolf's decomposition
of the pressure for the interior estimate and Liu-Wang's line for the boundary case.
%


At first,
let us introduce the definition of suitable weak solutions in the interior domain.
\begin{Definition} \label{Def:SWS}
Let  $\Om\subset\mathbb{R}^6$ be an open domain. $(u,\pi)$ is said to be a suitable weak solution to the steady Navier-Stoks equations (\ref{eq:SNS}) in $\Om$, if the following conditions hold.

(i)\,$u\in H^1(\Om),\,\pi\in L^{\frac 32}(\Om),\,\, f\in L^q(\Om)$,\,\, $q > 3$;

(ii)\,$(u,\pi)$ satisfies the equations (\ref{eq:SNS}) in the sense of distribution sense;

(iii)\,$u$ and $\pi$ satisfy the local energy inequality
\ben\label{eq:local energy}
2\int_{\Om}|\nabla u|^2\phi dx
\leq\int_{\Om}\big[|u|^2\triangle\phi+u\cdot\nabla\phi(|u|^2+2\pi)\big]+2 fu \phi dx,
\een
for any nonnegative $C^{\infty}$ test function $\phi$ vanishing at the boundary $\partial\Omega$ .
\end{Definition}

The existence of such a suitable weak solution can be found in \cite{FR3}. The major concern of this paper is the regularity and the main results can be stated as follows:

\begin{Theorem}\label{thm:c}
Let $(u,\pi)$ be a suitable weak solution to (\ref{eq:SNS}) in $B_1$.
Then $0$ is a regular point of $u$,
if there exists a small positive constant $\varepsilon$
such that the following conditions holds,
\beno
 r^{-3}\int_{B_r} |u(x)|^3 dx + r^{3q-6} \int_{B_r} |f(x)|^q dx <\varepsilon,
\eeno
for some $r\in (0,1)$.
\end{Theorem}

\begin{Remark}
The regularity criteria above for the 6D steady Navier-Stokes equations generalize recent interior regularity results by Dong-Strain \cite{DS}, where the pressure is small:
\beno
\int_{B_1} |u(x)|^3+|\pi(x)|^{\frac32} dx +  \int_{B_1} |f(x)|^2 dx\leq \varepsilon_0.
\eeno
%
\end{Remark}

Although the authors \cite{Wolf, CW, JWZ,WWZ} proved $\int_{Q_1}|u(x)|^3dx\leq \varepsilon_0$ in one scale can imply the regularity for the time-dependent Navier-Stokes equations, however it seems to be difficult for the regularity by only assuming  $\int_{Q_1}|\nabla u(x)|^2dx\leq \varepsilon_0$ in one scale.
Here for the steady equations, we have the following criterion:
\begin{Theorem}\label{thm:c'}
Let $(u,\pi)$ be a suitable weak solution to (\ref{eq:SNS}) in $B_1$.
Then $0$ is a regular point of $u$,
if there exists a small positive constant $\varepsilon$
such that the following conditions holds,
\beno
 \left(r^{-5}\int_{B_r} |u(x)| dx\right)^2 \left(r^{-2}\int_{B_r} |\nabla u(x)|^2 dx\right) + r^{-2}\int_{B_r} |\nabla u(x)|^2 dx+ r^{3q-6} \int_{B_r} |f(x)|^q dx <\varepsilon,
\eeno
for some $r\in (0,1)$.
\end{Theorem}

The theorem immediately implies the $2D$ Hausdorff measure of the set of singular points of $(u,\pi)$ in $B_1$ is equal to zero, and we omitted the proof, since it's standard as in \cite{DS}.


Second, let us introduce the definition of suitable weak solutions near the boundary.
\begin{Definition}
Let  $\Om\subset\mathbb{R}^6$ be an open domain, and $\Gamma\subset \partial\Omega$ be an open set. $(u,\pi)$ is said to be a suitable weak solution to the steady Navier-Stoks equations (\ref{eq:SNS}) in $\Omega$ near the boundary $\Gamma$, if the following conditions hold.

(i)\,$u\in H^1(\Omega),\,\pi\in L^{\frac 32}(\Omega),\, f\in L^{6}(\Omega)$;

(ii)\,$(u,\pi)$ satisfies the equations(\ref{eq:SNS}) in the sense of distribution sense and the boundary condition $u|_{\Gamma}=0$ holds;

(iii)\,$u$ and $\pi$ satisfy the local energy inequality
\ben\label{eq:local energy}
2\int_{\Omega}|\nabla u|^2\phi dx
\leq\int_{\Omega}\big[|u|^2\triangle\phi+u\cdot\nabla\phi(|u|^2+2\pi)\big]+2 fu \phi dx
\een
for any nonnegative $C^{\infty}$ test function $\phi$ vanishing at the boundary $\partial\Omega\backslash\Gamma$ .
\end{Definition}


Recall boundary regularity criteria in \cite{LW} stated as follows:
\begin{Proposition}[Theorem 1.2., Proposition 1.6., \cite{LW}]\label{prop:cd}
Let $(u,\pi)$ be a suitable weak solution to (\ref{eq:SNS}) in $B_1^+$ near the boundary $\{x\in B_1, x_6=0\}$. Then $0$ is a regular point of $u$, if there exists a small positive constant $\varepsilon_1$ such that one of the following conditions holds\\
(i)There exists ${\rho}_0>0$ such that
such that
\beno
\rho_0^{-3}\|u\|_{L^3(B_{\rho_0}^+)}^3+\rho_0^{-2}\|\nabla \pi\|_{L^{6/5}(B_{\rho_0}^+)}+\rho_0^3\|f\|_{L^3(B_{\rho_0}^+)}^3\leq\varepsilon_1,
\eeno
(ii)
\beno
\limsup_{r\rightarrow 0}r^{-2}\int_{B_r^+} |\nabla u(x)|^2 dx \leq \varepsilon_1,
\eeno
(iii)
\beno
\limsup_{r\rightarrow 0}r^{-3}\int_{B_r^+} | u(x)|^3 dx \leq \varepsilon_1.
\eeno
\end{Proposition}

The above result can be improved as follows:
\begin{Theorem}\label{thm:d}
Let $(u,\pi)$ be a suitable weak solution to (\ref{eq:SNS}) in $B_1^+$ near the boundary $\{x\in B_1, x_6=0\}$.
Then $0$ is a regular point of $u$,
if there exists a small positive constant $\varepsilon$
such that one of the following conditions holds, \\
$(i)$
\beno
r^{-3}\int_{B_r^+} |u(x)|^3 dx + r^3 \int_{B_r^+} |f(x)|^3 dx < \varepsilon,
\eeno
for some $r \in (0,1)$; \\
$(ii)$
\beno
r^{-2}\int_{B_r^+} |\nabla u(x)|^2 dx + r^3 \int_{B_r^+} |f(x)|^3 dx < \varepsilon,
\eeno
for some $r \in (0,1)$.
\end{Theorem}

\begin{Remark}
The regularity criteria above for the 6D steady Navier-Stokes equations improve recent interior regularity results in \cite{LW} by removing the condition of the pressure, which also improve the result of \cite{DG2}.
\end{Remark}

%

%



The rest of the paper is organized as follows. In Section 2, we introduce some notations, some technical lemmas and local energy estimates.
In Section 3 and 4, we prove Theorem \ref{thm:c} and Theorem \ref{thm:c'}, respectively.
Section 5 is devoted to the proof of Theorem \ref{thm:d}.
In Section 6, we show that any suitable weak solution to the steady Navier-Stokes equations is a local suitable weak solution.

Throughout this article, $C_0$ denotes an absolute constant independent of $u,\rho,r$
and may be different from line to line.

\setcounter{equation}{0}

\section{Notations and some technical lemmas}

Let $(u,\pi)$ be a solution to the steady Navier-Stokes equations (\ref{eq:SNS}). Set the following scaling:
\ben\label{eq:scaling}
u^{\lambda}(x)={\lambda}u(\lambda x),\quad \pi^{\lambda}(x)={\lambda}^2\pi(\lambda x),\quad f^{\lambda}(x)={\lambda}^3f(\lambda x),
\een
for any $\lambda> 0,$ then the family $(u^{\lambda},\pi^{\lambda})$ is also a solution of (\ref{eq:SNS}) with $f$ replaced by $f^{\lambda}(x)$.
Now define some quantities which are invariant under the scaling (\ref{eq:scaling}):

$$A(r)=r^{-4}\int_{B_r}|u(x)|^2dx,\quad C(r)=r^{-3}\int_{B_r}|u(x)|^3dx;$$

$$E(r)=r^{-2}\int_{B_r}|\nabla u(x)|^2dx;$$

$$D(r)=r^{-3}\int_{B_r}|\pi-(\pi)_{B_r}|^{\f32}dx,\quad (\pi)_{B_r}=\frac{1}{|B_r|}\int_{B_r}\pi dx;$$

$$F(r)=r^{3q-6}\int_{B_r}|f(x)|^qdx,$$
where
$B_r(x_0)$ is the  ball of radius $r$ centered at $x_0$, and
we denote $B_r(0)$ by $B_r$.
Moreover, a solution $u$ is said to be regular at $x_0$ if $u\in L^\infty(B_r(x_0))$ for some $r>0$.

Let us
introduce Wolf's pressure decomposition as in \cite{Wolf}. Given a bounded $C^2$-domain $G\subset R^n$ and $1<s<\infty$, define the operator
$
E_G:  W^{-1,s}(G)\rightarrow  W^{-1,s}(G)
$ as follows.
By the $L^p-$ theory of the steady Stokes system \cite{GSS}, for any $ F\in W^{-1,s}(G)$ there exists a unique pair $(v,\pi)\in W^{1,s}_{0}\times L^s_0(G)$ which solves the steady Navier-Stokes equations in the weak sense
\begin{equation}\label{eq:stokes}
\left\{\begin{array}{llll}
-\Delta v+\nabla \pi=F,\quad in \quad G\\
{\rm div }~ v=0, \quad in \quad G\\
v=0, \quad on \quad \partial G,
\end{array}\right.
\end{equation}
where $\pi\in L^s_0(G)$ denotes
\beno
\int_G \pi dx=0, \quad \pi\in L^s(\Omega).
\eeno
Then let $E_G(F)=\nabla \pi$, where $\nabla \pi$ denotes the gradient functional in $W^{-1,s}(G)$ defined by
\beno
<\nabla p,\psi>=-\int_Gp\nabla\cdot \psi dx,\quad \psi\in W^{1,s'}_0(G).
\eeno

The operator $E_G$ is bounded from $W^{-1,s}(G)$ into itself with $E_G(\nabla \pi)=\nabla \pi$ for all
$\pi\in L^s_0(G)$, and
\ben\label{eq:bound of wolf}
\|\pi\|_{L^s(G)}\leq C\|F\|_{W^{-1,s}(G)}.
\een
The norm of $E_G$ depends only on $s$ and the geometric properties of $G$, and is independent
of $G$, if $G$ is a ball or an annulus, which is due to the scaling properties of the Stokes equation.

Let us introduce the definition of local suitable weak solutions.
\begin{Definition} \label{Def:local SWS}
Let a bounded $C^2$-domain $\Om\subset\mathbb{R}^6$. $(u,\pi)$ is said to be a local suitable weak solution to the steady Navier-Stoks equations (\ref{eq:SNS}) in $\Om$, if the following conditions hold.

(i)\,$u\in H^1(\Om),\,\pi\in L^{\frac 32}(\Om),\,\, f\in L^{q}(\Om)$,\,\, $q > 3$;

(ii)\,$(u,\pi)$ satisfies the equations (\ref{eq:SNS}) in the sense of distribution sense;

(iii)\,for any ball $B\subset \Omega$, let $u$ and $\pi$ satisfy the local energy inequality
\ben\label{eq:local energy}
2\int_{\Om}|\nabla u|^2\phi dx
\leq\int_{\Om}\big[|u|^2\triangle\phi+u\cdot\nabla\phi(|u|^2+2\pi_{1}+2\pi_2)\big]+2 fu \phi dx
\een
for any nonnegative $C^{\infty}$ test function $\phi$ vanishing at the boundary $\partial B$, where
\beno
\nabla \pi_1=-E_B(u\cdot\nabla u),\quad \nabla \pi_2=E_B(\triangle u).
\eeno
\end{Definition}

\begin{Remark} \label{Rem:LSWS}
A suitable weak solution $(u,\pi)$ of (\ref{eq:SNS}) is a local suitable weak solution under the Definition \ref{Def:local SWS}. We prove this remark on Sec.6.
\end{Remark}


More precisely, we will prove the following proposition, which implies Theorem \ref{thm:c}.
\begin{Proposition}\label{prop: regularity}
Let $(u,\pi)$ be a local suitable weak
solution in $B_{1}$ to the Navier-Stokes equations (\ref{eq:SNS}). There exists absolute positive numbers $C_*$ and $\varepsilon$ such that  if
\beno
\int_{B_{1}}|u|^3dx + \left(\int_{B_{1}} |f|^q dx\right)^\frac3q \leq \varepsilon^3,
\eeno
then we have
\ben\label{eq:rk estimate}
r_k^{-6}\int_{B_{r_k}}|u|^3dz\leq C_*^3\varepsilon^3,\quad \forall~k\in N,
\een
where $r^k=2^{-k}$.
\end{Proposition}

Under the scaling (\ref{eq:scaling}), we also can define some quantities as follow:

$$A^+(r)=r^{-4}\int_{B_r^+}|u(x)|^2dx,\quad C^+(r)=r^{-3}\int_{B_r^+}|u(x)|^3dx;$$

$$E^+(r)=r^{-2}\int_{B_r^+}|\nabla u(x)|^2dx;$$

$$D^+(r)=r^{-3}\int_{B_r^+}|\pi-\pi_{B_r^+}|^{\f32}dx,\quad \pi_{B_r^+}=\frac{1}{|B_r^+|}\int_{B_r^+}\pi dx;$$

$$F^+(r)=r^{3}\int_{B_r^+}|f(x)|^3dx,$$

%
%
%

We need the following revised local energy inequality stated in \cite{LW}.
\begin{Proposition}\label{prop:local half space}
Let $0 < 16r < \rho \leq r_0$. It holds
\beno
&&k^{-2}A^+(r) + E^+(r)\\
&&\leq Ck^4\left(\frac r\rho\right)^2 A^+(\rho) + Ck^{-1} \left(\frac \rho r\right)^3 [C^+(\rho) + (C^+(\rho))^\frac 13 (D^+(\rho))^\frac23] \\
&&+ C \left(\frac\rho r\right)^2 (C^+(\rho))^\frac13 (F^+(\rho))^\frac13.
\eeno
Here $1 \leq k \leq \frac \rho r$ and constant $C$ is independent on $k,r,\rho$.
\end{Proposition}

\section{Interior regularity: proof of Theorem \ref{thm:c}}

In this section, we present the proof of Proposition \ref{prop: regularity}, whose proof is divided into several steps, which implies  Theorem \ref{thm:c}.
In details, we shall prove the key inequality (\ref{eq:rk estimate}) in Proposition \ref{prop: regularity} by using a strong induction
argument on $k$. Let $C_*$ be a constant which will be specified at the final moment.
From the definition of a local suitable weak solution the following local energy
inequality holds true for every nonnegative $\phi\in C_0^\infty(B_{\frac34})$,
\ben\label{eq:local energy2}
2\int_{B_{\frac34}}|\nabla u|^2\phi dx
\leq\int_{B_{\frac34}}\big[|u|^2\triangle\phi+u\cdot\nabla\phi(|u|^2+2\pi_{1}+2\pi_2)\big]+2 fu \phi dx.
\een

First, we introduce the following lemmas.

\begin{Lemma}[Cacciopolli type inequality]\label{lem:Cacciopolli}
Let $(u,\pi)$ be a local suitable weak
solution in $B_{1}$ to the Navier-Stokes equations (\ref{eq:SNS}). Then for any $0<R\leq 1$ there holds
\ben\label{ine:energy}
\|\nabla u\|_{L^2(B_{R/2})}^2\leq CR^{-2}\|u\|_{L^2(B_{R})}^2+CR^{-1}\|u\|_{L^3(B_{R})}^3+CR^{2}\|f\|_{L^2(B_{R})}^2.
\een
\end{Lemma}
\no
{\it Proof of Lemma \ref{lem:Cacciopolli}.}
For any $0 < R \leq \frac34$, choose $\phi = 1$ in $B_\tau$ and $\phi=0$ on $B_\rho^c$ with $\frac R2\leq \tau<\rho\leq R$ and
\beno
\nabla \pi_1=-E_{B_{\rho}}(u\cdot\nabla u),\quad \nabla \pi_2=E_{B_{\rho}}(\triangle u).
\eeno
It follows from (\ref{eq:local energy2}) and
 (\ref{eq:bound of wolf}) that
\beno
\int_{B_\tau} |\nabla u|^2 dx &\leq& C (\rho-\tau)^{-2} \int_{B_R} |u|^2 dx + C (\rho-\tau)^{-1} \int_{B_R}  |u|^3 dx
\\&&+ C (\rho-\tau)^{-1} \left(\int_{B_R} |u|^3 dx\right)^\frac 13 \left(\int_{B_\rho} |\pi_1|^\frac 32 dx\right)^\frac 23\\ &&+ C (\rho-\tau)^{-1} \left(\int_{B_R} |u|^2 dx\right)^\frac 12 \left(\int_{B_\rho} |\pi_2|^2 dx\right)^\frac 12+C\int_{B_R}|u||f|dx \\
&\leq& C (\rho-\tau)^{-2} \int_{B_R} |u|^2 dx + C (\rho-\tau)^{-1} \int_{B_R} |u|^3 dx \\
&&+ \frac 12 \int_{B_\rho} |\nabla u|^2 dx+ C\int_{B_R}|u||f|dx .
\eeno
By a standard iteration argument, the proof is complete.

Similar as Lemma 2.9 in \cite{CW} or Lemma 2.3 in \cite{JWZ}, we have the decay estimate of the pressure part $\pi_1$.
\begin{Lemma}[The pressure estimate]\label{lem:pressure1}
Let $(u,\pi)$ be a local suitable weak
solution in $B_{1}$ to the Navier-Stokes equations (\ref{eq:SNS}). Assume that for any $x_0\in B_{\frac12}$ and  $0< r \leq \frac12$ there holds
\beno
\int_{B_r(x_0)}|u\otimes u-({u\otimes u})_{B_r(x_0) }|^{\frac32}dx\leq CC_*^3 r^{6} \int_{B_1}|u|^3dx
\eeno
then for all $x_0 \in B_\frac12$ and $0 < r < \frac18$,
\beno
\int_{B_r(x_0)}|\pi_1-(\pi_1)_{B_r(x_0) }|^{\frac32}dx\leq CC_*^3 r^{6} \int_{B_1}|u|^3dx,
\eeno
where $\nabla \pi_1 = - E_{B_{\frac34}}(u \cdot \nabla u)$.
\end{Lemma}

\no
{\it Proof of Lemma \ref{lem:pressure1}.}
Without loss of generality, let $x_0 = 0$. Assume that $0 < \theta < \frac 14$ and $r \in (0,\frac18)$ fixed. Since $\nabla \pi_1 = - E_{B_{\frac34}}(u \cdot \nabla u)$, we can write
\beno
- \Delta v + \nabla (\pi_1 - (\pi_1)_{B_r}) = \nabla \cdot (u \otimes u) \quad {\rm in} \quad B_\frac34.
\eeno
Let
\beno
p_{0,r} = \frac{\partial_i \partial_j}{\Delta} ((u_i u_j - (u \otimes u)_{B_r})\zeta),
\eeno
where $\zeta = 1$ in $B_{r/2}$, $\zeta = 0$ on $B_r^c$, and
\beno
p_{h,r} = \pi_1 - (\pi_1)_{B_r} - p_{0,r}.
\eeno
Then $\pi_1 - (\pi_1)_{B_r} = p_{h,r} + p_{0,r}$ in $B_r$.
Using triangle inequality, we have
\beno
\int_{B_{\theta r}} |\pi_1 - (\pi_1)_{B_{\theta r}}|^\frac 32 dx &\leq& C \int_{B_{\theta r}} |p_{h,r} - (p_{h,r})_{B_{\theta r}}|^\frac 32 dx \\
&&+ C \int_{B_{\theta r}} |p_{0,r} - (p_{0,r})_{B_{\theta r}}|^\frac 32 dx \\
&& : = J_1 + J_2.
\eeno

For $J_1$, it follows from the properties of harmonic functions that
\beno
J_1 \leq C (\theta r)^\frac 32 \int_{B_{\theta r}} |\nabla p_{h,r}|^\frac 32 dx \leq C \theta^\frac {15}2 \int_{B_r} |p_{h,r}|^\frac 32 dx.
\eeno
For $J_2$, using Calder\'{o}n-Zygmund estimate, we have
\beno
J_2 &\leq& C \int_{B_{\theta r}} |p_{0,r}|^\frac 32 dx \leq C \int_{\mathbb{R}^6} |p_{0,r}|^\frac 32 dx \\
&\leq& C \int_{\mathbb{R}^6} |(u \otimes u - (u \otimes u)_{B_r}) \zeta|^\frac32 dx \\
&\leq& C \int_{B_r} |u \otimes u - (u \otimes u)_{B_r}|^\frac 32 dx.
\eeno

Combining these estimates, we get
\beno
\int_{B_{\theta r}} |\pi_1 - (\pi_1)_{B_{\theta r}}|^\frac 32 dx &\leq& C \theta^\frac{15}2 \int_{B_r} |p_{h,r}|^\frac 32 dx + C \int_{B_r} |u \otimes u - (u \otimes u)_{B_r}|^\frac 32 dx \\
&\leq& C \theta^\frac{15}2 \int_{B_r} |\pi_1 - (\pi_1)_{B_r}|^\frac 32 dx + C \int_{B_r} |u \otimes u - (u \otimes u)_{B_r}|^\frac 32 dx.
\eeno
Using standard iteration argument,
\ben\label{ine:lemma1}
\int_{B_r} |\pi_1 - (\pi_1)_{B_{r}}|^\frac 32 dx \leq C r^6 \int_{B_\frac 14} |\pi_1|^\frac 32 dx + CC_*^3 r^{6} \int_{B_1}|u|^3dx.
\een
Noting that $\nabla \pi_1 = -E_{B_\frac 34} (\nabla \cdot (u \otimes u))$, we have
\ben\label{ine:lemma2}
\int_{B_\frac 14} |\pi_1|^\frac 32 dx \leq C \int_{B_\frac 34} |u \otimes u - (u \otimes u)_{B_r}|^\frac 32 dx.
\een
Combining (\ref{ine:lemma1}) and (\ref{ine:lemma2}),  the proof is complete.
\endproof

\no
{\it Proof of Proposition \ref{prop: regularity}.} Let $r_n=2^{-n}$ and we
introduce a smooth function as
$$\Gamma_{n+1}(x)=\frac{1}{(r_{n+1}^2+|x-x_0|^2)^{2}},$$
which clearly satisfies
$$\triangle\Gamma_{n+1}=\frac{-24r_{n+1}^2}{(r_{n+1}^2+|x-x_0|^2)^4}<0.$$
Moreover, let
\beno
\chi(x)=1,\quad as \quad x\in B_{r_4}(x_0),
\eeno
and
\beno
\chi(x)=0,\quad as \quad x\in B_{r_3}^c(x_0).
\eeno

Obviously, the estimate of $(\ref{eq:rk estimate})$  holds for $k=1.$ Next
we assume that $(\ref{eq:rk estimate})$ holds for $k=1,\cdots,n.$

Taking the test function $\phi=\Gamma_{n+1}\chi$ in the local energy inequality (\ref{eq:local energy2}), we obtain that
\beno
&&-\int_{B_{r_3}(x_0)}|u|^2\chi\triangle\Gamma_{n+1} dx+2\int_{B_{r_3}(x_0)}|\nabla u|^2\chi\Gamma_{n+1} dx\\
&&\leq \int_{B_{r_3}(x_0)}|u|^2(\Gamma_{n+1}\triangle\chi+2\nabla\Gamma_{n+1}\cdot\nabla\chi)dx\\
&&+\int_{B_{r_3}(x_0)}u\cdot\nabla\phi |u|^2dx+2\int_{B_{r_3}(x_0)}u\cdot\nabla\phi ~\pi_1dx+2\int_{B_{r_3}(x_0)}u\cdot\nabla\phi ~\pi_2dx
\\
&&+2\int_{B_{r_3}(x_0)}f u \chi\Gamma_{n+1}  dx=I_1+\cdots+I_5.
\eeno
It follows from some straightforward computations that
\ben\label{eq:estimate of Gamma}
&&i)~\chi\Gamma_{n+1}(x,t)\geq C_0(r_{n+1})^{-4},\quad-\chi\triangle\Gamma_{n+1}(x,t)\geq C_0(r_{n+1})^{-6}\quad{\rm in} \,\, B_{r_{n+1}},\nonumber\\
&&ii)~|\nabla\phi|\leq |\nabla\Gamma_{n+1}|\chi+\Gamma_{n+1}|\nabla\chi|\leq C_0(r_{n+1})^{-5}\quad{\rm in} \,\, B_{\rho},\nonumber\\
&&iii)~|\Gamma_{n+1}\triangle\chi|+2|\nabla\Gamma_{n+1}\cdot\nabla\chi|\leq C_0\rho^{-6}\quad{\rm in} \,\, B_{\rho},
\een

{\bf \underline{Estimate of $I_1$}.} It follows from iii) of (\ref{eq:estimate of Gamma}) that
\beno
I_1\leq C_* \left(\int_{B_1}|u|^3dx\right)^{\frac23}.
\eeno

{\bf \underline{Estimate of $I_2$}.} Due to $|\nabla\phi|\leq Cr_k^{-5}$ in $B_{r_k}(x_0)\setminus B_{r_{k+1}}(x_0)$, we have
\beno
I_2=\int_{B_{r_3}(x_0)}u\cdot\nabla\phi |u|^2&\leq& \sum_{k=3}^n \int_{B_{r_k}(x_0)\setminus B_{r_{k+1}}(x_0) }|u|^3|\nabla\phi|+\int_{B_{r_{n+1}}(x_0)}|u|^3|\nabla\phi|\\
&\leq& C \sum_{k=3}^{n+1} r_k^{-5}\int_{B_{r_k}(x_0) }|u|^3dx\\
&\leq& CC_*^3\int_{B_1}|u|^3dx.
\eeno

{\bf \underline{Estimate of $I_3$}.} As in \cite{CKN}, we choose a series of cut-off functions $\chi_k$ satisfying
\begin{align*} \chi_k(x)=\left\{
\begin{aligned}
&1,\quad x\in  B_{r_{k+1}}(x_0),\\
&0, \quad x\in B_{r_k}(x_0)^c,
\end{aligned}
\right.
\end{align*}
for $k=3,\cdots,k+1.$
Then
\beno
\frac12I_3&=&\int_{B_{r_3}(x_0)}u\cdot\nabla\phi \pi_1dx\\
&\leq&\sum_{k=3}^n\int_{B_{r_k}(x_0)\setminus B_{r_{k+2}}(z_0)}(\pi_1-(\pi_1)_{B_{r_k}(x_0)})u\cdot \nabla[\phi(\chi_k-\chi_{k+1})]\\
&&+\int_{B_{r_2}(x_0)}(\pi_1)u\cdot \nabla[\phi(1-\chi_{3})]\\
&&+\int_{B_{r_{n+1}}(x_0)}(\pi_1-(\pi_1)_{B_{r_{n+1}}(x_0)})u\cdot \nabla[\phi\chi_{n+1}]=J_1+J_2+J_3.
\eeno
Due to $|\nabla (\phi(\chi_k - \chi_{k+1}))|\leq Cr_k^{-5}$ in $B_{r_k}(x_0)\setminus B_{r_{k+2}}(x_0)$, we have
\beno
J_1\leq C_* C\sum_{k=3}^n r_k^{-5} r_k^2\left(\int_{B_1}|u|^3dx\right)^{1/3}\|(\pi_1-(\pi_1)_{B_{r_k}(x_0)})\|_{L^{\frac32}(B_{r_k}(x_0))}.
\eeno
Since $\nabla\pi=E_{B_{\frac34}}(-u\cdot\nabla u)$ and
\beno
\int_{B_r(x_0)}|u\otimes u-({u\otimes u})_{B_r(x_0) }|^{\frac32}dx\leq CC_*^3 r^{6} \int_{B_1}|u|^3dx,
\eeno
then Lemma \ref{lem:pressure1} implies
\beno
\int_{B_r(x_0)}|\pi_1-(\pi_1)_{B_r(x_0) }|^{\frac32}dx\leq CC_*^3 r^{6} \int_{B_1}|u|^3dx,
\eeno
and
\beno
\|(\pi_1-(\pi_1)_{B_{r_k}(x_0)})\|_{L^{\frac32}(B_{r_k}(x_0))}\leq C r_k^{4}C_*^2\|u\|_{L^3(B_{\frac12})}^{2}.
\eeno

Hence we have
\beno
J_1\leq  CC_*^3\int_{B_1}|u|^3dx,
\eeno
and the other terms are similar.

{\bf \underline{Estimate of $I_4$}.} We still use the functions $\chi_k$.
\beno
I_4&=&\int_{B_{r_3}(x_0)}u\cdot\nabla\phi \pi_2dx\\
&\leq&\sum_{k=3}^n\int_{B_{r_k}(x_0)\setminus B_{r_{k+2}}(z_0)}(\pi_2-(\pi_2)_{B_{r_k}(x_0)})u\cdot \nabla[\phi(\chi_k-\chi_{k+1})]\\
&&+\int_{B_{r_2}(x_0)}\pi_2 u\cdot \nabla[\phi(1-\chi_{3})]\\
&&+\int_{B_{r_{n+1}}(x_0)}(\pi_2-{(\pi_2)}_{B_{r_{n+1}}(x_0)})u\cdot \nabla[\phi\chi_{n+1}]=J_1'+J_2'+J_3'.
\eeno
and by the induction assumption we get
\beno
J_1'\leq C_* C\sum_{k=3}^n r_k^{-5} r_k^2\left(\int_{B_1}|u|^3dx\right)^{1/3}r_k\|\pi_2-(\pi_2)_{B_{r_k}(x_0)}\|_{L^2(B_{r_k}(x_0))}.
\eeno
Due to the harmonic property of $\pi_2$, we have
\beno
\|\pi_2-{\pi_2}_{B_{r_k}(x_0)}\|_{L^2(B_{r_k}(x_0))}\leq Cr_k^4\|\pi_2\|_{L^2(B_{\frac12})}\leq Cr_k^4[\|u\|_{L^3(B_{\frac34})}^{\frac32}+\|u\|_{L^3(B_{\frac34})} + ||f||_{L^2(B_\frac34)}].
\eeno
where we used the local energy inequality.
And the other terms are similar.

Hence, we have
\beno
I_4\leq CC_*[\|u\|_{L^3(B_{\frac34})}^{\frac52} + \|u\|_{L^3(B_{\frac34})}^2 + \|u\|_{L^3(B_{\frac34})}||f||_{L^q(B_\frac34)}].
\eeno

{\bf \underline{Estimate of $I_5$}.}
Since $|(\chi_k - \chi_{k+1}) \Gamma_{n+1}| \leq C r_k^{-4}$ in $B_{r_k}(x_0)\setminus B_{r_{k+2}}(x_0)$, we have
\beno
\frac12 I_5&=&\int_{B_{r_3}(x_0)}f u \chi\Gamma_{n+1}  dx\\
&\leq& \sum_{k=3}^n \int_{B_{r_k}(x_0)\setminus B_{r_{k+2}}(x_0)}f u \chi(\chi_k-\chi_{k+1})\Gamma_{n+1}  dx+\int_{B_{r_3}(x_0)}f u \chi(1-\chi_{3})\Gamma_{n+1}\\
&&+\int_{B_{r_3}(x_0)}f u \chi(\chi_{n+1})\Gamma_{n+1}=J_1''+J_2''+J_3''.
\eeno
Since $q > 3$, we have
\beno
J_1''\leq C C_*\sum_{k=3}^n r_k^{-4}r_k^{6-\frac 6q} \|u\|_{L^3(B_1)} \|f\|_{L^q(B_1)}\leq  C C_* \|u\|_{L^3(B_1)}\|f\|_{L^q(B_1)}.
\eeno
The estimate of other term is similar as $J_1''$. Hence we have
\beno
I_5 \leq C C_* \|u\|_{L^3(B_1)}\|f\|_{L^q(B_1)}.
\eeno

Combining $I_1$, $\dots$, $I_5$, we have
\beno
&&r_{n+1}^{-6}\int_{B_{r_{n+1}}(x_0)}|u|^2dx+r_{n+1}^{-4}\int_{B_{r_{n+1}}(x_0)}|\nabla u|^2dx\\
&\leq& CC_*[C_*^2\|u\|_{L^3(B_{\frac12})}^{3}+\|u\|_{L^3(B_{\frac12})}^2]+C C_* \|u\|_{L^3(B_1)}\|f\|_{L^q(B_1)},
\eeno
which implies that
\beno
r_{n+1}^{-6}\int_{B_{r_{n+1}}(x_0)}|u|^3dx\leq  CC_*^{\frac32}[C_*^3\|u\|_{L^{3}(B_{1})}^{\frac92}+\|u\|_{L^3(B_{1})}^3]+C C_*^{\frac32} \|u\|_{L^3(B_1)}^{\frac32}\|f\|_{L^q(B_1)}^{\frac32}.
\eeno
Then by choosing $\frac 52 C \leq C_*^\frac32$ and $\varepsilon$ small such that $C_*^{3}\|u\|_{L^3(B_1)}^\frac32 \leq \frac12$, we have
\beno
r_{n+1}^{-6}\int_{B_{r_{n+1}}(x_0)}|u|^3dx\leq  C_*^3\left(\int_{B_1}|u|^3dx+\|f\|_{L^q(B_1)}^3\right).
\eeno

The proof is complete.\endproof

\section{Proof of Theorem \ref{thm:c'}}

{\it Proof of Theorem \ref{thm:c'}.}
By Sobolev's embedding theorem, for $0<r<\rho$ we have
\beno
r^{-3}\int_{B_r} |u|^3 dx &\leq& Cr^{-3}\int_{B_r} |u - (u)_{B_\rho}|^3 dx + Cr^{-3} \int_{B_r} |(u)_{B_\rho}|^3 dx \\
&\leq& C r^{-3}\left(\int_{B_\rho} |\nabla u|^2 dx\right)^\frac32 + C r^3 \rho^{-18} \left(\int_{B_\rho} |u| dx\right)^3\\
&\leq& C \left(\frac \rho r\right)^3 E(\rho)^\frac32 + C \left(\frac r \rho\right)^3 \left(\rho^{-5} \int_{B_\rho} |u| dx\right)^3.
\eeno

{\bf Case I:} If $\rho^{-5} \int_{B_\rho} |u| dx \leq E^\frac12(\rho)$, we have
\beno
r^{-3}\int_{B_r} |u|^3 dx \leq C \big[\left(\frac \rho r\right)^3 + \left(\frac r \rho\right)^3\big] E^\frac32(\rho).
\eeno
Choosing $r = \frac 12 \rho$, noting that the assumption of Theorem \ref{thm:c'}, we have $C(r) \leq \varepsilon$.

{\bf Case II:} If $\rho^{-5} \int_{B_\rho} |u| dx > E^\frac12(\rho)$, let $r = \theta \rho$, we have
\beno
C(\theta \rho) \leq C \theta^{-3} E(\rho)^\frac 32 + C \theta^3 \left(\rho^{-5} \int_{B_\rho} |u| dx\right)^3.
\eeno
Choosing $\theta^6 = \frac{E(\rho)^\frac32}{\left(\rho^{-5} \int_{B_\rho} |u| dx\right)^3}$, we have $\theta < 1$ and
\beno
E(\theta \rho) \leq C E(\rho)^\frac 34 \left(\rho^{-5} \int_{B_\rho} |u| dx\right)^\frac32.
\eeno

Applying Theorem \ref{thm:c}, the proof of Theorem \ref{thm:c'} is complete.
\endproof

\section{Boundary regularity: proof of Theorem \ref{thm:d}}

In this section, we follow the same line as in \cite{LW} to prove the boundary regularity, and one major difference is that we only need to assume that the smallness condition holds on a fixed ball without the pressure term. Our new observation is based on the suitable pressure decomposition and the new estimates of the Stokes system including global and interior estimates (for example, see  \cite{Tsai-2018}, \cite{Ka}).

First of all, let us recall global Stokes estimates with zero boundary condition.
\begin{Lemma}[Theorem 2.13, \cite{Tsai-2018}]\label{lem:global Stokes estimate}
Let $\Omega$ be a bounded domain in $\mathbb{R}^n$, $n \geq 2$ and $q \in (1,\infty)$. For every $f=(f_{ij}) \in L^q(\Omega)$, there is a unique $q-$weak solution $v \in W_0^{1,q}(\Omega)$ of
\ben\label{eq:Stokes1}
- \Delta v_j + \nabla_j p = \partial_i(f_{ij}); \quad {\rm div}~v = 0,
\een
satisfying
\beno
||v||_{W^{1,q}(\Omega)} +||p-(p)_{\Omega}||_{L^{q}(\Omega)}\leq C ||f||_{L^q(\Omega)},
\eeno
where the constant $C$ only depends on $q$ and $\Omega$.
For every $g \in L^q(\Omega)$, there is a unique $q-$weak solution $v \in W^{2,q} \cap W_0^{1,q}(\Omega)$ of
\ben\label{eq:Stokes2}
- \Delta v + \nabla p = g; \quad {\rm div}~v = 0,
\een
satisfying
\beno
||v||_{W^{2,q}(\Omega)} +||\nabla p||_{L^{q}(\Omega)}\leq C ||g||_{L^q(\Omega)},
\eeno
where the constant $C$ only depends on $q$ and $\Omega$.
\end{Lemma}

\begin{Remark}
A vector field $v$ defined on $\Omega$ is a very weak solution of (\ref{eq:Stokes1}), if $v \in L_{loc}^2(\Omega)$ and $v$ satisfies
\beno
\int_{\Omega} \nabla v : \nabla \zeta =- \langle f,\nabla\zeta \rangle; \quad \int_{\Omega} v \cdot \nabla \phi = 0;
\eeno
for all $\zeta \in C_c^\infty(\Omega)$ with $\nabla \cdot \zeta = 0$ and for all $\phi \in C_c^\infty(\Omega)$.
A very weak solution $v$ is a $q-$weak solution of (\ref{eq:Stokes1}), if $v \in W^{1,q}(\Omega)$. The $q-$weak solution of (\ref{eq:Stokes2}) is similar.
\end{Remark}

The interior estimate for the pressure plays an important role in the following arguments. At this time, one main feature is that the velocity is zero only on part of the boundary, but the estimation of the higher derivative of pressure can not depend on the lower derivative of pressure. We recall a theorem by Kang  in \cite{Ka} as follows.
\begin{Lemma}[Theorem 3.8, \cite{Ka}]\label{lem:pressure near boundary}
Let $\Omega \subset \mathbb{R}^n$ be a domain of class $\mathcal{C}^{k+2}$ and $k$ be an integer with $-1 \leq k < \infty$ and $1 < q <\infty$. Suppose $g \in W^{k,q}(\Omega_{r_0})$ and $u \in W^{1,q}(\Omega_{r_0})$ with a unique pressure satisfying $ \int_{\Omega_{r_0} } p=0$ solve the following Stokes system:
\begin{align*}
\left\{
\begin{aligned}
- \Delta u + \nabla p = g,& \quad {\rm in} ~~\Omega_{r_0} \\
\nabla \cdot u = 0,& \quad {\rm in} ~~\Omega_{r_0} \\
u = 0,& \quad {\rm on} ~~B_{r_0} \cap \partial\Omega,\\
\end{aligned}
\right.
\end{align*}
in a weak sense. Let $r,s$ be positive numbers with $0 \leq r < s \leq r_0$. Then the following estimate holds:
\beno
||u||_{W^{k+2,q}(\Omega_r)} + ||p||_{W^{k+1,q}(\Omega_r)} \leq C \left(||g||_{W^{k,q}(\Omega_{r_0})} + ||u||_{L^1(\Omega_s)}\right),
\eeno
where $C = C(k,n,q,r,s,\Omega)$ and $\Omega_r = \Omega \cap B_r$ with $r \leq r_0$. Here $r_0$  is comparable to the radius of the sphere contained within this domain $\Omega$.
\end{Lemma}


Next we prove Theorem \ref{thm:d}.

\no
{\it Proof of Theorem \ref{thm:d}.} The proof of this theorem is divided into three parts.

{\bf Step I: The pressure estimate.}
First, we choose a domain $\tilde{B}^+$ with a smooth boundary such that $B_{\frac 34}^+ \subset \tilde{B}^+ \subset B_1^+$.
 Let $\tilde{B}^+_\rho = \{\rho x: x \in \tilde{B}^+\}$, which implies $\tilde{B}^+_\rho$ is also smooth.
For $0<\rho<1$, let $v$ and $\pi_1$ be the unique solution to the following boundary value problem of Stokes system
\begin{eqnarray}\nonumber
 \left\{
    \begin{array}{llll}
    \displaystyle - \Delta v + \nabla \pi_1 = f - u \cdot \nabla u& \quad {\rm in} \quad \tilde{B}^+_\rho,\\
    \displaystyle {\rm{div}} ~ v = 0& \quad {\rm in} \quad \tilde{B}^+_\rho,\\
    \displaystyle v = 0& \quad {\rm on} \quad \partial\tilde{B}^+_\rho,\\
    \displaystyle (\pi_1)_{\tilde{B}^+_\rho} = \int_{\tilde{B}^+_\rho} \pi_1 dx = 0.\\
    \end{array}
 \right.
\end{eqnarray}
Due to the uniqueness of the linear Stokes system, we decompose $v = v_1 + v_2$ and $\pi_1 = \pi_{11} + \pi_{12}$ in this way, which satisfy
\begin{eqnarray}\nonumber
 \left\{
    \begin{array}{llll}
    \displaystyle - \Delta v_1 + \nabla \pi_{11} = f& \quad {\rm in} \quad \tilde{B}^+_\rho,\\
    \displaystyle {\rm{div}} ~ v_1 = 0& \quad {\rm in} \quad \tilde{B}^+_\rho,\\
    \displaystyle v_1 = 0& \quad {\rm on} \quad \partial\tilde{B}^+_\rho,\\
    \displaystyle (\pi_{11})_{\tilde{B}^+_\rho} = \int_{\tilde{B}^+_\rho} \pi_{11} dx = 0,\\
    \end{array}
 \right.
\end{eqnarray}
and
\begin{eqnarray}\nonumber
 \left\{
    \begin{array}{llll}
    \displaystyle - \Delta v_2 + \nabla \pi_{12} = - u \cdot \nabla u& \quad {\rm in} \quad \tilde{B}^+_\rho,\\
    \displaystyle {\rm{div}} ~ v_2 = 0& \quad {\rm in} \quad \tilde{B}^+_\rho,\\
    \displaystyle v_2 = 0& \quad {\rm on} \quad \partial\tilde{B}^+_\rho,\\
    \displaystyle (\pi_{12})_{\tilde{B}^+_\rho} = \int_{\tilde{B}^+_\rho} \pi_{12} dx = 0.\\
    \end{array}
 \right.
\end{eqnarray}
With the help of Lemma \ref{lem:global Stokes estimate}, we get
\ben\label{ine:estimate of pi_1} \nonumber
&\rho^{-2}||v_1||_{L^{\frac32}(\tilde{B}^+_\rho) }+ \rho^{-1}||\pi_{11}||_{L^{\frac32}(\tilde{B}^+_\rho)} \leq C ||f||_{L^{\frac32}(\tilde{B}^+_\rho)},\\
&\rho^{-1}||v_2||_{L^{\frac32}(\tilde{B}^+_\rho)} + ||\pi_{12}||_{L^{\frac32}(\tilde{B}^+_\rho)} \leq C |||u|^2||_{L^{\frac32}(\tilde{B}^+_\rho)},
\een
where the constant $C$ is independent of $\rho$ due to the scaling transform. Then it follows that
\ben\label{ine:pi 1 estimate}\nonumber
||\pi_1||_{L^\frac32(\tilde{B}^+_\rho)} &\leq& ||\pi_{11}||_{L^\frac32(\tilde{B}^+_\rho)} + ||\pi_{12}||_{L^\frac32(\tilde{B}^+_\rho)} \\
&\leq& C \rho^3 ||f||_{L^3(B^+_\rho)} + C ||u||_{L^3(B^+_\rho)}^2,
\een
and
\beno
||v||_{L^\frac32(\tilde{B}^+_\rho)} &\leq& ||v_1||_{L^\frac32(\tilde{B}^+_\rho)} + ||v_2||_{L^\frac32(\tilde{B}^+_\rho)} \\
&\leq& C \rho^4 ||f||_{L^3(\tilde{B}^+_\rho)} + C \rho ||u||_{L^3(\tilde{B}^+_\rho)}^2.
\eeno

%

On the other hand, let $w = u - v$, $\pi_2 = \pi - (\pi)_{\tilde{B}^+_\rho} - \pi_1$, then $\int_{\tilde{B}^+_\rho} \pi_2 = 0$. Moreover, ($w$, $\pi_2$) solves the following boundary value problem:
\begin{eqnarray}\nonumber
 \left\{
    \begin{array}{llll}
    \displaystyle - \Delta w + \nabla \pi_2 = 0& \quad {\rm in} \quad \tilde{B}^+_\rho,\\
    \displaystyle {\rm{div}} ~ w = 0& \quad {\rm in} \quad \tilde{B}^+_\rho,\\
    \displaystyle w = 0& \quad {\rm on} \quad \partial \tilde{B}^+_\rho \cap \{x_6 = 0\}.\\
    \end{array}
 \right.
\end{eqnarray}
Using Lemma \ref{lem:pressure near boundary} by choosing $r_0=\rho$, $r=\frac14 \rho$ and $s=\frac12\rho$, we have
\ben\label{ine:pi 2 estimate}\nonumber
\rho^{3-\frac6q}||\nabla \pi_2||_{L^q(B_{\frac14\rho}^+)} &\leq& C \rho^{-5}||w||_{L^1(B_{\frac12\rho}^+)} \leq C \rho^{-5} ||u||_{L^1(B_{\frac12\rho}^+)} + C \rho^{-5}||v||_{L^1(B_{\frac12\rho}^+)} \\
&\leq& C \rho^{-1} ||u||_{L^3(B_{\frac12\rho}^+)} + C \rho ||f||_{L^3(\tilde{B}^+_\rho)} + C \rho^{-2} ||u||_{L^3(\tilde{B}^+_\rho)}^2 ,
\een
where the constant $C$ is independent of the radius $\rho$.

Combining (\ref{ine:pi 1 estimate}) and (\ref{ine:pi 2 estimate}), for $0<4r<\rho$ we have
\beno
||\pi - (\pi)_{B^+_r}||_{L^\frac32(B_r^+)} &\leq& ||\pi_1 - (\pi_1)_{B^+_r}||_{L^\frac32(B_r^+)} + ||\pi_2 - (\pi_2)_{B^+_r}||_{L^\frac32(B_r^+)} \\
&\leq& C ||\pi_1||_{L^\frac32(\tilde{B}^+_\rho)} + Cr^{5-\frac6q}||\nabla \pi_2||_{L^q(B_r^+)} \\
&\leq& C \rho^3 ||f||_{L^3(\tilde{B}^+_\rho)} + C ||u||_{L^3(\tilde{B}^+_\rho)}^2 + C \left(\frac r\rho\right)^{5-\frac6q} \rho ||u||_{L^3(B_{\frac12\rho}^+)} \\
&&+ C \left(\frac r\rho\right)^{5-\frac6q} \rho^3 ||f||_{L^3(\tilde{B}^+_\rho)} + C \left(\frac r\rho\right)^{5-\frac6q} ||u||_{L^3(\tilde{B}^+_\rho)}^2,
\eeno
which yields that for $0<4r<\rho<1$ there holds
\ben\label{ine:D + r}
D^+(r) &\leq& C \left(\frac r\rho\right)^{\frac92-\frac9q} (C^+(\rho))^\frac12 + C \left(\frac \rho r\right)^3 (F^+(\rho))^\frac12 + C \left(\frac \rho r\right)^3 C^+(\rho)
\een

{\bf Step II:  Theorem \ref{thm:d} under the assumption $(i)$.}
Proposition \ref{prop:local half space} and (\ref{ine:D + r}) tell us that for all $0<\theta<\frac14$, the following estimate holds:
\beno
k^{-2} A^+(\theta \rho) + E^+(\theta \rho) &\leq& C_0 k^4 \theta^2 A^+(\rho) + C_0 k^{-1} \theta^{-3} C^+(\rho) + C_0 k^{-1} \theta^{-3} (C^+(\rho))^\frac13 (D^+(\rho))^\frac23 \\
&&+ C_0 \theta^{-2} (C^+(\rho))^\frac13 (F^+(\rho))^\frac13,
\eeno
and
\ben\label{ine:D + theta rho}\nonumber
D^+(\theta \rho) &\leq& C_0 \theta^{\frac92-\frac9q} (C^+(\rho))^\frac12 + C_0 \theta^{-3} (F^+(\rho))^\frac12 + C_0 \theta^{-3} C^+(\rho)
\een
where $k \in [1,\theta^{-1}]$ and $C_0$ is a constant independent of $\rho$ and $\theta$.


Let $G(r) = k^{-2} A^+(r) + E^+(r) + \gamma^{-1} (D^+(r))^\frac43$, where $\gamma > 0$ to be decided.  By the embedding inequality
\ben\label{eq:embedding}
C^+(\rho) \leq C_0 (E^+(\rho))^\frac32,
\een
since $u=0$ on partial boundaries, we have
\beno
G(\theta\rho) &\leq& C_0 k^4 \theta^2 A^+(\rho) + C_0 k^{-1} \theta^{-3} C^+(\rho) + C_0 k^{-1} \theta^{-3} (C^+(\rho))^\frac13 (D^+(\rho))^\frac23 \\
&&+ C_0 \theta^{-2} (C^+(\rho))^\frac13 (F^+(\rho))^\frac13 + C_0 \gamma^{-1} \theta^{6-\frac{12}q} E^+(\rho) \\
&&+ C_0 \gamma^{-1} \theta^{-4} (F^+(\rho))^\frac23 + C_0 \gamma^{-1} \theta^{-4} (C^+(\rho))^\frac43 \\
&\leq&  C_0 k^4 \theta^2 A^+(\rho) +  C_0 k^{-1} \theta^{-3} C^+(\rho) +\frac14 G(\rho) + C_0 \gamma k^{-2} \theta^{-6} (C^+(\rho))^{\frac23} \\
&&+ C_0 \theta^{-4} (F^+(\rho))^\frac23 + C_0 \gamma^{-1} \theta^{6-\frac{12}q} E^+(\rho) + C_0 \gamma^{-1} \theta^{-4} (F^+(\rho))^\frac23 \\
&&+ C_0 \gamma^{-1} \theta^{-4} (C^+(\rho))^{\frac43} \\
&\leq& \frac14 G(\rho) + C_0 \left(k^6\theta^2  + \gamma k^{-2} \theta^{-6} + \gamma^{-1} \theta^{6-\frac{12}q}\right) G(\rho) \\
&&+ C_0 k^{-1} \theta^{-3} C^+(\rho) + C_0 \theta^{-4} (F^+(\rho))^\frac23 +  C_0 \gamma^{-1} \theta^{-4} (F^+(\rho))^\frac23 \\
&&+ C_0 \gamma^{-1} \theta^{-4} (C^+(\rho))^{\frac43}
\eeno
Letting $k = \theta^{-\frac14}$ and $\gamma = \theta^{6-\frac14}$,  we have
\beno
k^6\theta^2 + \gamma k^{-2} \theta^{-6} + \gamma^{-1} \theta^{6-\frac{12}q} \leq \theta^\frac12 + \theta^\frac14 + \theta^{\frac14-\frac{12}q}.
\eeno
Take $q > 48$ and $\theta = \theta_0$, which  satisfy
\ben\label{ine:theta 0}
C_0\left(\theta_0^\frac12 + \theta_0^\frac14 + \theta_0^{\frac14-\frac{12}q}\right) < \frac14,
\een
which implies
\ben\label{eq:G estimate}
G(\theta_0 \rho) \leq \frac12 G(\rho) + C_1\left(C^+(\rho) + (C^+(\rho))^{\frac43}  + (F^+(\rho))^\frac23\right),
\een
where the constant $C_1$ is only dependent on $C_0$ and $\theta_0$.

Without loss of generality, under the condition of $(i)$, assume that $C^+(\rho_0) + F^+(\rho_0) < \varepsilon_0<1$ with $0<\rho_0<1$. Then (\ref{ine:D + r})
implies that $D^+(\theta_0^{\frac12}\rho_0) \leq   C_1(\varepsilon_0^\frac12  + \varepsilon_0)$ with $C_1$ only depending on $\theta_0$. On the other hand,
by Proposition \ref{prop:local half space}
\beno
&&A^+(r) + E^+(r) \leq C_0 \frac{r^2}{\rho^2} A^+(\rho) \\
&&+ C_0 \left(\frac\rho r\right)^3 \left[C^+(\rho) + (C^+(\rho))^\frac13  (D^+(\rho))^\frac23 \right]+ C_0 \left(\frac\rho r\right)^2 (C^+(\rho))^\frac13 (F^+(\rho))^\frac13,
\eeno
which, the H\"{o}lder inequality of $A^+(\rho)\leq C_0(C^+(\rho))^\frac23 $ and the estimate of $D^+(\theta_0^{\frac12}\rho_0)$ imply that
\ben\label{eq:G1}
G(\rho_1)&=&G(\theta_0\rho)=\theta_0^{\frac12} A^+(\theta_0\rho) + E^+(\theta_0\rho) + \theta_0^{-\frac{23}{4}} (D^+(\theta_0\rho))^\frac43\nonumber\\
&\leq& C_1(\varepsilon_0^{\frac23}+ \varepsilon_0^{\frac43} +\varepsilon_0)\leq\varepsilon_1
\een
where  $\rho_1 = \theta_0 \rho_0$,  $\theta_0$ is decided by (\ref{ine:theta 0}) and $C_1$ only depending on $\theta_0$.

It follows from (\ref{eq:G estimate}) and (\ref{eq:G1}) that
\beno
G(\theta_0 \rho_1) \leq \frac12G(\rho_1) + C_1(\varepsilon_1^2 + \varepsilon_1^\frac32 + \varepsilon_0^\frac23).
\eeno
where we used (\ref{eq:embedding}) and the monotonicity of $F^+(\rho)$ with respect to $\rho$. Here $C_1$ only depends on $\theta_0$.
Choosing $\varepsilon_1 = \varepsilon_1(\varepsilon_0)$ small enough, we have
\beno
G(\theta_0\rho_1) \leq \frac12\varepsilon_1 +  C_1(\varepsilon_1^2 + \varepsilon_1^\frac32 + \varepsilon_1) \leq \varepsilon_1^\frac12.
\eeno
Assume that $G(\theta_0^j \rho_1) \leq \varepsilon_1^\frac12$ holds for $1\leq j\leq k$, next we verify the case of $k+1.$ Then there holds
\beno
G(\theta_0^{k+1} \rho_1) &\leq& \frac12G(\theta_0^{k}\rho_1) + C_1(\varepsilon_1 + \varepsilon_1^\frac34 + \varepsilon_0^\frac23)\leq \varepsilon_1^\frac12
\eeno
where we used (\ref{eq:G estimate}), (\ref{eq:embedding}) and the monotonicity of $F^+(\rho)$ again.
Hence by mathematical induction,  for all $j \in \mathbb{N}$ there holds
\ben\label{ine:induction}
G(\theta_0^j \rho_1) \leq \varepsilon_1^\frac12.
\een
Consequently,
for any $r \in (0,\theta_0\rho_1)$, there exists a constant $j$ such that $\theta_0^{j} \rho_1 \leq r < \theta_0^{j-1} \rho_1$. Thus
\beno
E^+(r) \leq  \theta_0^{-2} E^+(\theta_0^{j-1} \rho_1)\leq \theta_0^{-2} G(\theta_0^{j-1} \rho_1)
\eeno
Noting that (\ref{ine:induction}), we have
\beno
E^+(r) \leq C_3 \varepsilon_1^\frac12, \quad \forall ~~r \in (0,\theta_0\rho_1),
\eeno
where $C_3$ only dependent on $\theta_0$ is a absolute constant. Using Proposition 1.6, the proof of Theorem \ref{thm:d} is complete under the assumption $(i)$.

{\bf Step III:  Theorem \ref{thm:d} under the assumption $(ii)$.}

For the assumption $(ii)$, it's obvious from (\ref{eq:embedding}) that
\beno
C^+(\rho) \leq C_0 (E^+(\rho))^\frac32,
\eeno
which implies the assumption $(i)$.
\endproof


\section{Appendix: the proof of  Remark \ref{Rem:LSWS}}

We clarify that any suitable weak solution to the steady Navier-Stokes equations is a local suitable weak solution.

\begin{Lemma} \label{lem:LSWS SWS}
Let $(u,\pi)$ as in Definition \ref{Def:SWS} be a suitable weak solution to the Navier-Stokes equations (\ref{eq:SNS}). Then $u$ is a local suitable weak solution in the sense of Definition \ref{Def:local SWS}.
\end{Lemma}

\no {\bf Proof of Lemma \ref{lem:LSWS SWS}.}
Let $B \subset \mathbb{R}^6$ be a fixed ball. Without loss of generality, we assume that $f = 0$. Define
\beno
\nabla \pi_{0,B} = E_B(-u \cdot \nabla u + \Delta u).
\eeno
Since $E_B$ is a bounded operator and $\nabla \pi \in W^{-1,q}(B)$, we have
\beno
\nabla \pi &=& E_B(\nabla \pi) = E_B(-u \cdot \nabla u + \Delta u) = \nabla \pi_{0,B} \\
&=& E_B(-u \cdot \nabla u) + E_B(\Delta u) \\
&:=& \nabla \pi_1 + \nabla \pi_2.
\eeno
Since $(u,\pi)$ is a suitable weak solution, we have
\beno
2\int_{\Omega} |\nabla u|^2 \phi dx \leq \int_{\Omega} \big[ |u|^2 \triangle \phi + u \cdot \nabla \phi (|u|^2 + 2\pi) \big] dx.
\eeno
Applying integration by parts, it follows that
\beno
2\int_{\Omega} |\nabla u|^2 \phi dx &\leq& \int_{\Omega} \big[ |u|^2 \triangle \phi + u \cdot \nabla \phi (|u|^2 + 2\pi_{0,B}) \big] dx \\
&\leq& \int_{\Omega} \big[ |u|^2 \triangle \phi + u \cdot \nabla \phi (|u|^2 + 2\pi_1 + 2\pi_2) \big] dx.
\eeno
Thus the proof is complete.
\endproof

\bigskip

\noindent {\bf Acknowledgments.}
W. Wang was supported by NSFC under grant 11671067.


\end{document}